\documentclass[a4paper]{amsart}
\usepackage{amsfonts}
\usepackage{amssymb}
\usepackage{epsfig}
\usepackage{latexsym}
\usepackage{graphicx}

\theoremstyle{plain}
 \newtheorem{Th}{Theorem}
 \newtheorem{PROP}[Th]{Proposition}
 \newtheorem{LEMMA}[Th]{Lemma}

\theoremstyle{definition}
 \newtheorem{DEF}[Th]{Definition}

\theoremstyle{remark}
 
 \newtheorem{REM}[Th]{Remark}

\def\C{\mathcal C}
\def\M{\mathcal M}

\def\S{\mathcal S}
\def\T{\mathcal T}
\def\R{\mathcal R}

\let\longto\longrightarrow

\makeatletter
\def\@oddfoot{\footnotesize\normalfont \hfil-- \thepage\ --\hfil}
\def\@evenfoot{\footnotesize\normalfont \hfil-- \thepage\ --\hfil}
\makeatother

\makeatletter
  \textfloatsep\floatsep
  \intextsep\floatsep
\makeatother

\begin{document}

\title[THE COMPLEX OF PANT DECOMPOSITIONS OF A SURFACE]
      {THE COMPLEX OF PANT DECOMPOSITIONS OF A SURFACE}

\author[Silvia Benvenuti \and Riccardo Piergallini]
       {Silvia Benvenuti \and Riccardo Piergallini}
                        \address[Silvia Benvenuti
                                 \and Riccardo Piergallini]
                        {Dipartimento di Matematica e Informatica
                        -- Universit\`a di Camerino\\
                        Via Madonna delle Carceri 9\\
                        62032 Camerino -- Italia}
                        \email[S.~Benvenuti]
                              {silvia.benvenuti@unicam.it}
                        \email[R.~Piergallini]
                              {riccardo.piergallini@unicam.it}

\subjclass{Primary 57M50; Secondary 57M20, 30F60}
\keywords{Pant decomposition complex, mapping class group, complex of curves.}

\begin{abstract}
We exhibit a set of edges ({\em moves}) and 2-cells ({\em relations})
making the complex of pant decompositions on a surface a simply
connected complex. Our construction, unlike the previous ones, keeps
the arguments concerning the structural transformations independent from
those deriving from the action of the mapping class group. The moves and
the relations turn out to be supported in subsurfaces with $3g-3+n=1,2$
(where $g$ is the genus and $n$ is the number of boundary components),
illustrating in this way the so called Grothendieck principle.
\end{abstract}

\maketitle

\thispagestyle{empty}

\vspace{-.5cm}
\tableofcontents
\vspace{-.5cm}

\section{Introduction}
\label{introduction}

Let $\Sigma = \Sigma_{g,n}$ be a connected, compact, oriented surface
of genus $g$ with $n$ boundary components ($g,n \geq 0$). In order to
describe an algebraic or geometric object $\tau(\Sigma)$, it is often
convenient to represent $\Sigma$ as the result of gluing together
several simple pieces, which should be surfaces with boundary.

This happens for example in the study of the mapping class group
$\M(\Sigma) = \M_{g,n}$ (whose presentation may be obtained starting
from those of the mapping class group of some simple subsurfaces, as
proved in~\cite{Benv}), in the pantwise construction of hyperbolic
structures (Fenchel-Nielsen construction) and in the construction of
modular functors (defined by gluing the vector spaces associated to
simpler subsurfaces, provided they satisfy some {\em gluing axiom}).

Depending on the situation, it is convenient to choose the building
blocks for our surfaces from different {\em Lego boxes}: one may
choose, for example, a big Lego box (the {\em grande boite}
in~\cite{Groth}), whose pieces are all the spheres with any number
of boundary components, or maybe a smaller Lego box, containing only
spheres with at most three boundary components. In many cases, a
{\em cheap} Lego box, made of identical pieces (namely spheres with
three boundary components or hexagons) may be sufficient.

\begin{figure}[htb]
\epsfig{file=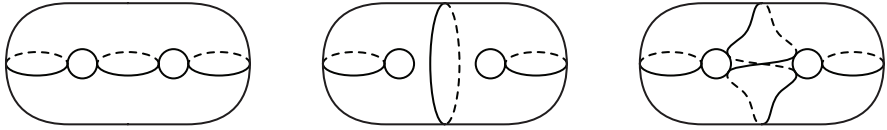}
\caption{Three different pant decompositions of $\Sigma_{2,0}$.}
\label{decompositions}
\end{figure}

No matter which box one chooses, it is evident that each surface
admits an infinite number of different decompositions with pieces
out of that box (see for instance Figure~\ref{decompositions}).

Thus, if one wants to describe an object $\tau(\Sigma)$ using a
decomposition $d$ of $\Sigma$ (i.e. computing $\tau(\Sigma)$ as a
$\tau(\Sigma,d)$), in order for this object to be well defined it 
is necessary to construct canonical isomorphisms between the objects
computed starting from different decompositions, i.e. to construct
isomorphisms $$f:\tau(\Sigma,d_1) \to \tau (\Sigma,d_2),$$ where $d_1$
and $d_2$ are any two different decompositions of $\Sigma$.

For instance, coming back to the cases mentioned before: while
studying the mapping class group it turns out that different
``slicing'' of $\Sigma$ produce different presentations for
$\M(\Sigma)$, and we look for a procedure to get any presentation
from any other; while constructing hyperbolic structures, different
pant decompositions lead to the description of different charts of
the atlas of the Teichm\"uller space of $\Sigma$, and we look for
the change of chart; while building a modular functor, different
ways of sewing $\Sigma$ result in different bases for the vector
space associated to $\Sigma$, and we need to write down the matrices
giving the change of basis ({\em duality matrices} in~\cite{MS}).

Therefore, once the Lego box is fixed, we have to describe the set of
all the decompositions of a surface into pieces from that Lego box,
considered up to isotopy, and the set of all the transformations between
different (non-isotopic) decompositions. More precisely, our aim is to
exhibit:

\begin{description}\rightskip\leftmargin\advance\rightskip\itemindent
 \item[elementary moves,] such that we can go from a given decomposition
   to any other through a sequence of these moves;
 \item[defining relations,] describing when a sequence of elementary
   moves applied to a decomposition yields the same decomposition.
\end{description}

Following the philosophy introduced by Hatcher and Thurston in their
pioneering paper \cite{Hatcher-Thurston}, such problem can be reformulated 
as follows: 
we consider all the decompositions of the surface $\Sigma = \Sigma_{g,n}$ 
up to isotopy, as the vertices of a 2-dimensional CW complex $\R(\Sigma) = 
\R_{g,n}$; then we put an edge between two vertices if the corresponding 
decompositions are related by one of our candidate {\em moves}, and we 
cup off a loop with a 2-cell if the corresponding sequence of moves is 
one of our candidate {\em relations}. Hence we are reduced to check that 
$\R_{g,n}$ is simply connected. Indeed, this complex is connected if and 
only if the set of our candidates moves is complete. Moreover, it is 
simply connected if and only if any relation between elementary moves
follows from the ones we have cupped off.
Actually in \cite{Hatcher-Thurston} the focus is on the {\em cut system
complex} $\C(\Sigma)$, but in the appendix the authors suggest that the 
same program could be carried over to the case of pant decompositions.

The {\em pant decomposition complex} $\R_{g,n}$ was studied in \cite{MS} 
by Moore and Seiberg.
Unfortunately, their proof of the connectedness and simply connectedness
of the so built complex contains some serious gaps. In particular, it is
based on the knowledge of an explicit presentation for the mapping class
groups $\M_{g,n}$, which was then unknown. Indeed, at the moment they
were writing (1989), the only known finite presentations were those of
the modular groups $\M_{0,n}$ and $\M_{g,0}$.

The mapping class group of the surface we are examinating enters the
playground since the elements of this group act as transformations
on the set of the decompositions of $\Sigma$: for instance, the
decomposition shown in the right hand picture of
Figure~\ref{decompositions} is obtained from the one in the center by
a Dehn twist along the dotted curve. Anyway, not all transformations
between decompositions are elements of the mapping class group: for
instance, the decomposition shown in the left hand picture of
Figure~\ref{decompositions} cannot be transformed into the center one
(nor into the right one) by any homeomorphism of $\Sigma$.
Therefore, the set of transformations between different decompositions 
of the surface $\Sigma$ contains a core, that is the mapping class group 
$\M(\Sigma)$, and something additional: the idea in~\cite{MS} is to get 
rid of this ``extra part'' and eventually come to the study of 
$\M(\Sigma)$, which is what Moore and Seiberg could not carry out.

More recently, this problem has been overcome by using the Cerf 
theoretic techniques introduced in \cite{Hatcher-Thurston}, either 
directly \cite{FG, HLS} or passing through a projection on the cut 
systems complex \cite{BK}.

In this paper, we come back to the original Moore and Seiberg's
approach and fill in the gaps, exploiting the presentations of the
mapping class groups that we have at present.

Indeed, in~\cite{Gervaisart}, Gervais provides a presentation of
$\M_{g,n}$, in terms of Dehn twists. Another presentation, as quotients
of Artin groups, is described by Matsumoto in~\cite{Matsumoto} for
$\M_{g,1}$ and then generalized to the case of $\M_{g,n}$ by Labru\`ere
and Paris in~\cite{luis-cat}.

A general machinery for getting presentations of the mapping class
groups in any preferred ``style'' (for example in terms of Dehn
twists or as quotients of Artin groups) is given in~\cite{Benv}. The
procedure introduced in that paper takes as input the well known
presentations for the {\em sporadic surfaces} ($\Sigma_{0,4}$,
$\Sigma_{1,1}$, $\Sigma_{0,5}$ and $\Sigma_{1,2}$) according to some
``style'' and returns a presentation in the same ``style'' for every
$\Sigma_{g,n}$. Such flexibility makes this last approach suitable
for obtaining different descriptions of the complex $\R_{g,n}$ we
may need in different contexts.

\medskip

The paper is organized as follows. Section~\ref{maintools} contains 
the first definitions and the main tools.
The construction of the complex $\R_{g,n}$ is subdivided into two
independent steps. Given a surface $\Sigma_{g,n}$, in
Section~\ref{CS} we consider a 2-dimensional, finite, simply
connected CW complex $\S_{g,n}$, whose vertices are in one-to-one
correspondence with the {\em combinatorial structures} of pant
decompositions of $\Sigma_{g,n}$, i.e. with the
$\M_{g,n}$-equivalence classes of such decompositions. Then, in
Section~\ref{decomposizioni in pantaloni}, by using a presentation
of $\M_{g,n}$, we construct an infinite simply connected complex
codifying all the pant decompositions on $\Sigma_{g,n}$ and the
transformations between them, that is the desired $\R_{g,n}$.

\medskip

When we started thinking of this Lego-Teichm\"uller game, our aim
was to understand something new about a Grothendieck conjecture on
the subject (see~\cite{Groth}). In our context this conjecture can
be expressed, roughly speaking, as follows. Take the family of all
the complexes $\R_{g,n}$ and stack them in levels, putting at level
$k$ all the $\R_{g,n}$ with $3g-3+n=k$. Then, to describe the whole
tower of complexes, it is sufficient to describe its first and
second floor, i.e. the complexes with $3g-3+n = 1,2$, that are
precisely the ones of the sporadic surfaces $\Sigma_{0,4}$,
$\Sigma_{1,1}$, $\Sigma_{0,5}$ and $\Sigma_{1,2}$.

It is worth noticing that we obtained this Grothendieck principle for
the mapping class groups as a byproduct, in~\cite{Benv}. Namely, we
proved that the generators and the relations which are needed to present
the mapping class group of any surface are supported in subsurfaces
living at the first and second Grothendieck floor. Combining this result
with the new ones described in the present paper, we are now able to
prove the Grothendieck conjecture in the above stated form.

\section{First definitions and main tools}
\label{maintools}

Let $\Sigma = \Sigma_{g,n}$ be a connected, compact, oriented surface, 
of genus $g$ with $n$ boundary components. A {\em pant decomposition} of
$\Sigma$ is a decomposition of the surface into a finite number of pants,
determined by a collection of disjoint simple closed curves in the interior
of $\Sigma$. We recall that a {\em pant} is a closed disk with two smaller
open disks removed, i.e. the surface $\Sigma_{0,3}$. As usual, the family
of curves and therefore the induced pant decomposition are always
considered up to isotopy.

\begin{figure}[htb]
\epsfig{file=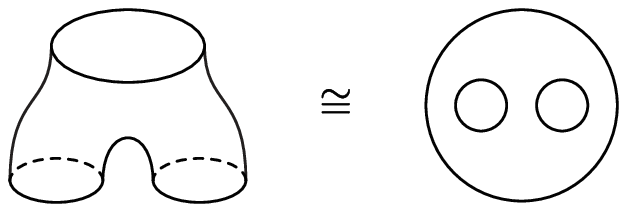}
\caption{A pant.}
\label{pant}
\end{figure}

More precisely, let $\alpha=\{\alpha_1,\dots,\alpha_k\}$ be a collection 
of pairwise disjoint closed loops on $\Sigma$. We denote by $\Sigma_\alpha$ 
the natural compactification of $\Sigma - \cup\{\alpha_i\}_{i=1,\dots,k}$\break
(obtained compactifying each component by the addiction of three boundary\break
curves), and by $\rho_\alpha:\Sigma_\alpha \to \Sigma$ the continuous map 
induced by the inclusion of $\Sigma -\cup\{\alpha_i\}_{i=1,\dots,k}$ into 
$\Sigma$.

\begin{DEF}\label{pant_decomposition}
We say that the family $\alpha=\{\alpha_1,\dots,\alpha_k\}$ determines a
{\em pant decomposition} of $\Sigma$ if each component of $\Sigma_\alpha$
is a pant, or equivalently if $\Sigma -\cup\{\alpha_i\}_{i=1,\dots,k}$
consists of $h$ components, each of which is homeomorphic to the
interior of a pant.
\end{DEF}

\begin{figure}[htb]
\epsfig{file=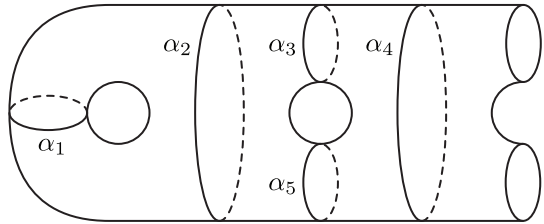}
\caption{A pant decomposition.}
\label{decomp_example}
\end{figure}

It can be easily shown that $\Sigma_{g,n}$ admits a pant decomposition
provided $$(g,n) \not\in \{(0,0),(0,1),(0,2),(1,0)\}.$$ Moreover, the
integers $k$ and $h$ are uniquely determined, and they are given by
$$k=3g-3+n \hspace{.5cm}\text{and}\hspace{0.5cm} h=2g-2+n=-\chi(\Sigma),$$
where $\chi(\Sigma)$ is the Euler characteristic of the surface $\Sigma$.

Let $N$ be a connected component of $\Sigma_\alpha$ (notations as above). 
We say that a boundary curve $\gamma$ of $N$ is {\em an exterior boundary 
curve} if $\rho_\alpha(\gamma)$ is a boundary component of $\Sigma$. 
For each curve $\alpha_i$ in the family $\alpha$ there are two distinct 
boundary curves $\gamma,\gamma'$ in $\Sigma_\alpha$ such that $\rho_\alpha(\gamma) = \rho_\alpha(\gamma') = \alpha_i$, and two possibilities arise:
either $\gamma$ and $\gamma'$ are boundary curves of the same connected
component $N$ of $\Sigma_\alpha$ (like $\alpha_1$ in 
Figure~\ref{decomp_example}), or $\gamma$ is a boundary component of $N$ 
and $\gamma'$ is a boundary component of a different connected component
$N'$ (like $\alpha_2$ in Figure~\ref{decomp_example}). In the first case 
we call $\alpha_i$ a {\em non-separating limit curve} of $N$, while in
the second case we call it a {\em separating limit curve} of $N$ and $N'$.

We denote by $\M_{g,n}$ the mapping class group of $\Sigma$, i.e. the
group of the isotopy classes of orientation preserving homeomorphisms
$h: \Sigma \to \Sigma$ which fix pointwise the boundary of $\Sigma$.
Clearly, we have an induced action of $\M_{g,n}$ on the set of the pant
decompositions of $\Sigma$. We call $\M_{g,n}$-equivalent two pant
decompositions that can be obtained from one another by such action.

\begin{DEF}
A {\em combinatorial structure} of pant decompositions of $\Sigma_{g,n}$
is a class of $\M_{g,n}$-equivalence of pant decompositions on 
$\Sigma_{g,n}$.
\end{DEF}

Now we can state the following proposition, whose proof is trivial.

\begin{PROP}\label{Gorbite} Let $[\alpha] = [\alpha_1, \ldots,
\alpha_{3g-3+n}]$ and $[\beta] = [\beta_1, \ldots, \beta_{3g-3+n}]$ be
two isotopy classes of curves defining pant decompositions of $\Sigma$.
Then the two pant decompositions belong to the same combinatorial
structure if and only if there exists a one-to-one correspondence
between the components of $\Sigma_\alpha$ and those of $\Sigma_{\beta}$
and there exists a permutation $\sigma \in S_{3g-3+n}$ such that, for
every pair $(N,N')$ where $N$ is any component of $\Sigma_\alpha$ and
$N'$ the corresponding component of $\Sigma_{\beta}$, we have:
\begin{itemize}\leftskip-15pt
\item[1)] if $\gamma$ is an exterior boundary curve of $N$ there exist an 
  exterior boundary curve $\gamma'$ of $N'$ such that $\rho_\alpha(\gamma)
  = \rho_{\beta }(\gamma')$ (i.e. $N$ and $N'$ have the same boundary
  components);
\item[2)] if $\alpha_i$ is a separating (resp. non-separating) limit
  curve of $N$, then $\beta_{\sigma(i)}$ is a separating (resp.
  non-separating) limit curve of $N'$.
\end{itemize}
\end{PROP}

In the light of the previous proposition, once a numbering of the
boundary components of $\Sigma_{g,n}$ is fixed, any combinatorial
structure of pant decompositions of $\Sigma_{g,n}$ can be encoded
by its dual graph. This graph has $2g-2+n$ trivalent vertices
corresponding to pants, and $n$ univalent vertices corresponding
to the boundary components. Moreover, the univalent vertices are
labelled by $\{1,2,\dots,n\}$ according to the fixed numbering
of the boundary components of $\Sigma_{g,n}$.

\medskip

The next proposition will be used as a criterion for the simply
connectedness of the complexes we are going to construct.

\begin{PROP}\label{teoremazero}
Let $\pi:C \to D$ be a surjective cellular map between 2-dimensional
CW complexes. Suppose the following conditions are satisfied:

\begin{itemize}\leftskip-15pt
\item[1)] for any vertex $v \in D$, the fiber $\pi^{-1}(v)$ is
  connected and simply connected in $C$;
\item[2)] for any oriented edge $e:v_1 \longto v_2$ in $D$ and any two
  liftings $e':v'_1 \longto v'_2$ and $e'':v''_1 \longto v''_2$ of $e$
  in $C$, there exist two paths $\gamma_1:v'_1 \longto v''_1$ in the fiber
  $\pi^{-1}(v_1)$ and $\gamma_2:v'_2 \longto v''_2$, $\gamma_i$ in the
  fiber $\pi^{-1}(C_2)$, such that the square
  \begin{center}
  \strut\kern-5mm\epsfig{file=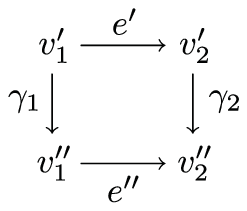}
  \end{center}
  is contractible in $C$.
\end{itemize}
Then, if $D$ is simply connected, $C$ is simply connected as well.
\end{PROP}

The proof of this proposition is straightforward, and it is left to
the reader.

\section{The combinatorial structures of pant decompositions}
\label{CS}

This section is devoted to the construction of a simply connected
complex $\S_{g,n}$ whose vertices represent the combinatorial
structures of pant decompositions of $\Sigma_{g,n}$ or equivalently
their dual graphs, as we said in the previous section. Moreover, in
the last subsection we will lift the complex $\S_{g,n}$ to another
one, denoted by $\widetilde{\S}_{g,n}$, codifying all of {\em decorated combinatorial structures}, i.e. the combinatorial structures of pant
decompositions whose curves are ordered. This is a technical tool which 
will be needed in Section~\ref{decomposizioni in pantaloni}.

We define $V(\S_{g,n})$ to be the set of all connected graphs with
$n$ univalent vertices, also called {\em free ends}, and $2g-2+n$
trivalent vertices. A standard computation shows that each of these
graphs has $3g-3+2n$ edges, $n$ connecting a free end to a trivalent
vertex, and the remaining $3g-3+n$ connecting two (possibly
coinciding) trivalent vertices.

\begin{figure}[htb]
\begin{center}
\epsfig{file=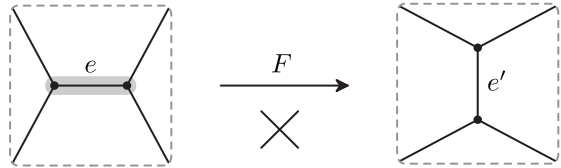}
\end{center}
\caption{The combinatorial $F$ move.}
\label{Fmove}
\end{figure}

On such graphs we consider the local move shown in Figure~\ref{Fmove}, 
that we call {\em combinatorial $F$ move}, according to the literature, 
as it can be thought as the fusion of two adjacent trivalent vertices 
followed by the inverse of a similar fusion. We warn the reader that 
the graphs in this picture, as well as in the following ones, should 
be considered as abstract graphs, regardless of their planar 
representation.

The picture means that the graph is unchanged outside a regular
neighborhood of an edge $e$ between two distinct trivalent vertices,
while such edge is replaced with a new edge $e'$, connecting two new
trivalent vertices. To be more precise, let $v_1$ and $v_2$ be the
trivalent vertices connected by the edge $e$ on which we perform the
$F$ move. Then, for each $i$, there are two (possibly coinciding)
edges other than $e$ having $v_i$ as a vertex. We label those
connected to $v_1$ by $a$ and $b$, and those connected to $v_2$ by
$c$ and $d$. Thus, the starting graph represents the coupling
$(ab)(cd)$. Therefore, for each edge between trivalent vertices
there are exactly two possible ways to perform the $F$ move (i.e.
two possible results for the $F$ move), corresponding to the two
possible changes of coupling:
$$(ab)(cd) \to (ad)(bc);$$
$$(ab)(cd) \to (ac)(bd).$$

\smallskip

The $F$ move is oriented and in the pictures we mark the starting edge 
$e$ on which the move is performed by surrounding it by a grey region.

To emphasize that the $F$ move can be realized as the contraction of
the edge $e$ followed by the inverse of a similar contraction of the
edge $e'$, we will label the arrows representing $F$ moves by the
corresponding intermediate graphs.

\begin{figure}[b]
\begin{center}
\epsfig{file=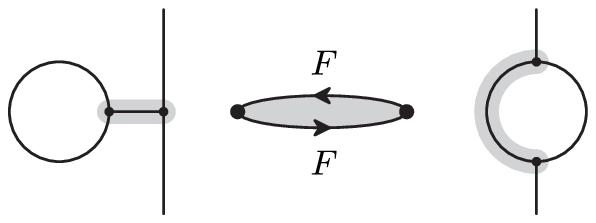}
\end{center}
\caption{A bigon.} 
\label{bigon}
\end{figure}

\begin{figure}[htb]
\begin{center}
\epsfig{file=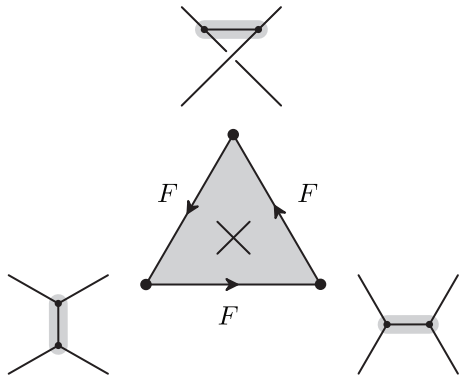}
\end{center}\vskip-3pt
\caption{A triangle.}
\label{triangle}
\end{figure}

\begin{figure}[htb]
\begin{center}\vskip3pt
\epsfig{file=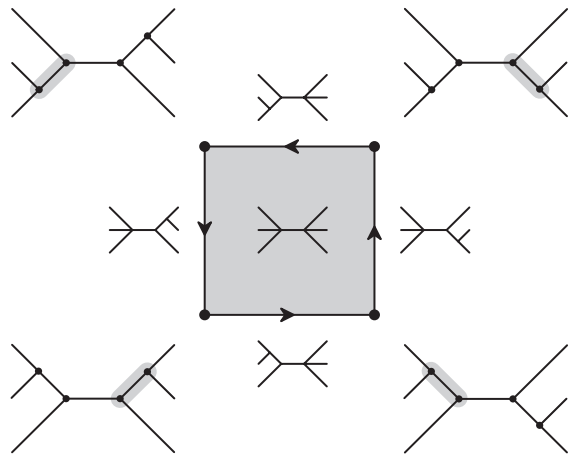}
\end{center}\vskip-3pt
\caption{A DC square.}
\label{square}
\end{figure}

\begin{figure}[b]
\begin{center}\vskip-18pt
\epsfig{file=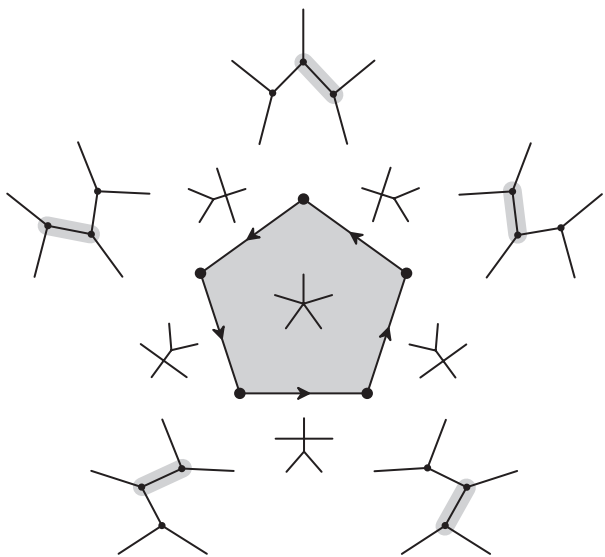}
\end{center}\vskip-3pt
\caption{A pentagon.}
\label{pentagon}
\end{figure}

Now, we define $\S_{g,n}$ to be the complex having $V(\S_{g,n})$ as
the set of vertices, an undirected edge connecting any two vertices
which are related by a combinatorial $F$ move, and the following
2-cells:

\break

\begin{description}\rightskip\leftmargin\advance\rightskip\itemindent
\item[bigons,] as in Figure~\ref{bigon};
\item[triangles,] as in Figure~\ref{triangle};
\item[squares] of ``disjoint commutativity'' (DC squares),
  as in Figure~\ref{square};
\item[pentagons,] as in Figure~\ref{pentagon}.
\end{description}

\smallskip

Coherently with what we have done for the edges, we consider only
one two cell for each loop of edges as in the Figures, even if the
same loop may be represented by different sequences of $F$ moves.

\begin{REM}
 Notice that triangles appear for
$2g-2+n \geq 2$, squares and pentagons show up when $2g-2+n \geq 3$,
while bigons are required for $g \geq 1$.
\end{REM}

\begin{REM}
The arrows appearing in the pictures do not represent an orientation
for the edges, but they are only intended to specify the oriented
$F$ move we are considering.
\end{REM}

The main result of this section is the following.

\begin{Th}\label{MMM}
The complex $\S_{g,n}$, with $2g-2+n \geq 1$, is simply connected.
\end{Th}

The proof takes the next two subsections and proceeds by an
induction scheme based on the diagram below, where the complexes
$\S_{g,n}$ are staked in layers corresponding to the value of
$2g-2+n$.

\smallskip

\begin{center}
\strut\kern4mm\epsfig{file=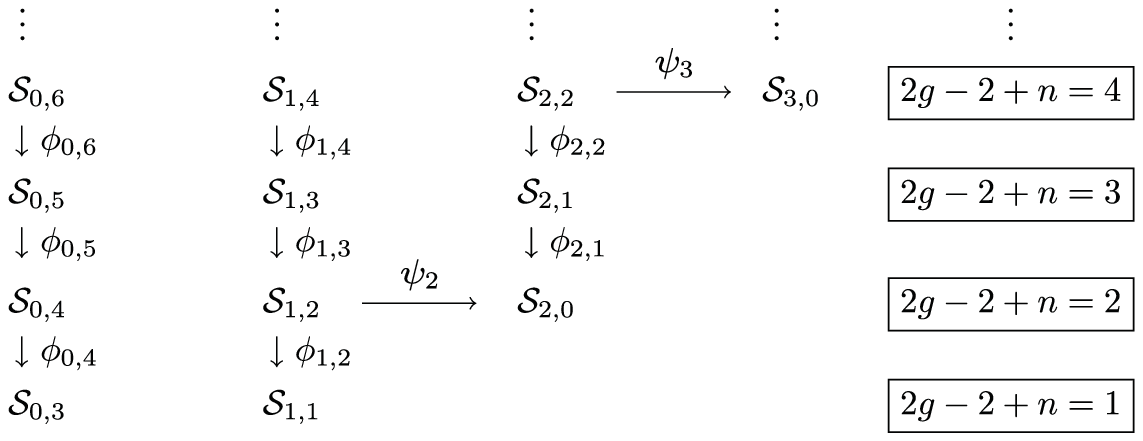}
\end{center}

\medskip

The base for the inductive argument is provided by the two complexes
$\S_{0,3}$ and $\S_{1,1}$, which both consist of a single vertex and
are therefore trivially connected and simply connected.

In Subsection~\ref{riduzione_componenti_di_bordo} we define the maps
$\phi_{g,n}:\S_{g,n} \to \S_{g,n-1}$ and show that they satisfy the
hypotheses of Proposition \ref{teoremazero}, in order to derive the
simply connectedness of $\S_{g,n}$ from that of $\S_{g,n-1}$
(Proposition~\ref{RIDUZIONE}). On the other hand, the maps
$\psi_g:\S_{g-1,2} \to \S_{g,0}$ are defined in
Subsection~\ref{casochiuso} and are used to prove that the simply
connectedness of $\S_{g-1,2}$ implies that of $\S_{g,0}$
(Proposition~\ref{CASO_CHIUSO}).

\subsection{Reducing the number of boundary components}
\label{riduzione_componenti_di_bordo}\strut\nobreak\medskip

We start by defining the map $$\phi_{g,n}:\S_{g,n} \to \S_{g,n-1}\;.$$

\begin{figure}[htb]
\begin{center}
\epsfig{file=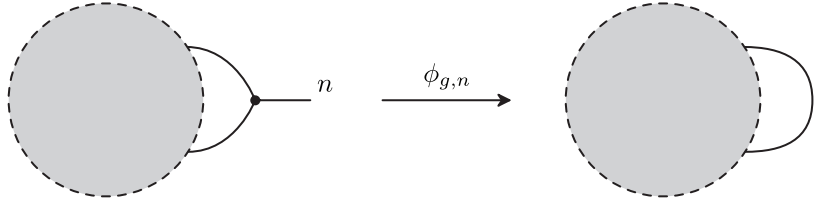}
\end{center}
\caption{Definition of $\phi_{g,n}$ on $V(\S_{g,n})$.}
\label{phi_vert}
\end{figure}

If $\Gamma$ is a vertex of $\S_{g,n}$, $\phi_{g,n}(\Gamma)$ is the graph
obtained from $\Gamma$ by contracting the last free end as in
Figure~\ref{phi_vert}: the result is a connected graph with $n-1$ free 
ends in the set $\{1,\dots,n-1\}$ and $2g-2+n-1$ trivalent vertices, 
i.e. a vertex of $\S_{g,n-1}$.

\begin{figure}[htb]
\begin{center}
\epsfig{file=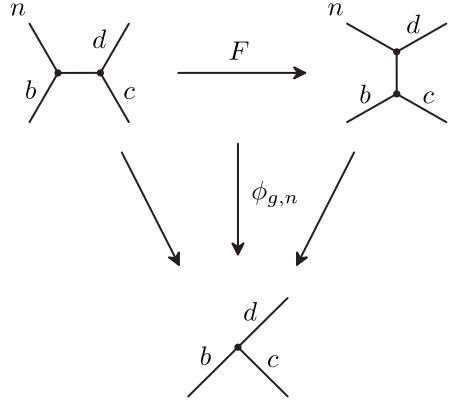}
\end{center}
\caption{Defining $\phi_{g,n}$ on the edges of $\S_{g,n}$.}
\label{phi_edge}
\end{figure}

As far as the edges are concerned, with the notations introduced at
the beginning of this section, we have two possibilities: either $n
\not \in \{a,b,c,d\}$ or $n\in \{a,b,c,d\}$. In the first case,
$\phi_{g,n}(F)$ is defined as the $F$ move between
$\phi_{g,n}(\Gamma_1)$ and $\phi_{g,n}(\Gamma_2)$. In the second
case, being $\phi_{g,n}(\Gamma_1)=\phi_{g,n}(\Gamma_2)$, we can
define $\phi_{g,n}(F)=\phi_{g,n}(\Gamma_i)$, as depicted in
Figure~\ref{phi_edge}.

Finally we define $\phi_{g,n}$ on the 2-cells in the following way.
If $\phi_{g,n}$ sends the boundary of a $2$-cell of $\S_{g,n}$ onto
the boundary of a $2$-cell of $\S_{g,n-1}$ having the same shape,
then $\phi_{g,n}$ sends the $2$-cell of $\S_{g,n}$ in the
corresponding one of $\S_{g,n-1}$. Otherwise, it may happen that an
edge of the $2$-cell collapses to a vertex of $\S_{g,n-1}$, as in
Figure~\ref{phi_edge}. If this happens for a bigon or a triangle,
then $\phi_{g,n}$ sends the entire boundary to the same vertex, and
we send the 2-cell itself to such vertex. On the other hand, if it
happens for a square or a pentagon, then the image of the boundary
reduces to an edge, and we send all the 2-cell to such an edge.

The map $\phi_{g,n}$ defined in this way is obviously cellular and
surjective both on vertices and on edges. To see that it is also
surjective on $2$-cells, we observe that each $2$-cell of
$\S_{g,n-1}$ is the image of a $2$-cell of $\S_{g,n}$ of the same
shape. In fact, the $F$ moves on the boundary of a $2$-cell of
$\S_{g,n-1}$ always leave a free edge on which the new trivalent
vertex can be created, to get the boundary of a $2$-cell of
$\S_{g,n}$ such that no edge collapses.

Now, we use $\phi_{g,n}$ to get the following inductive step for the
proof of Theorem \ref{MMM}.

\begin{PROP}\label{RIDUZIONE}
If the complex $\S_{g,n-1}$ is simply connected, then also $\S_{g,n}$
is simply connected, for every $g$ and $n$ such that $2g-2+n\geq 2$.
\end{PROP}

\begin{proof}
The thesis follows by applying Proposition~\ref{teoremazero} to the map 
$\phi_{g,n}$, once we prove that the required conditions concerning the
preimages of vertices and edges are fulfilled.

\smallskip
\noindent {\bf\em Claim 1.} {\em The fiber over each vertex is connected
and simply connected in $\S_{g,n}$.}
\smallskip

Let $\Gamma$ be a vertex of $\S_{g,n-1}$. Then the vertices of
$\phi_{g,n}^{-1}(\Gamma)$ consists of all the graphs with $n$ free ends
obtained from $\Gamma$ by inserting an edge with one free end labelled by
$n$ and the other end creating a new trivalent vertex along any of the
$3g-3+2(n-1)$ edges of $\Gamma$, which is therefore split into two edges
(see Figure~\ref{phi_inv_vert}).

\begin{figure}[htb]
\begin{center}
\epsfig{file=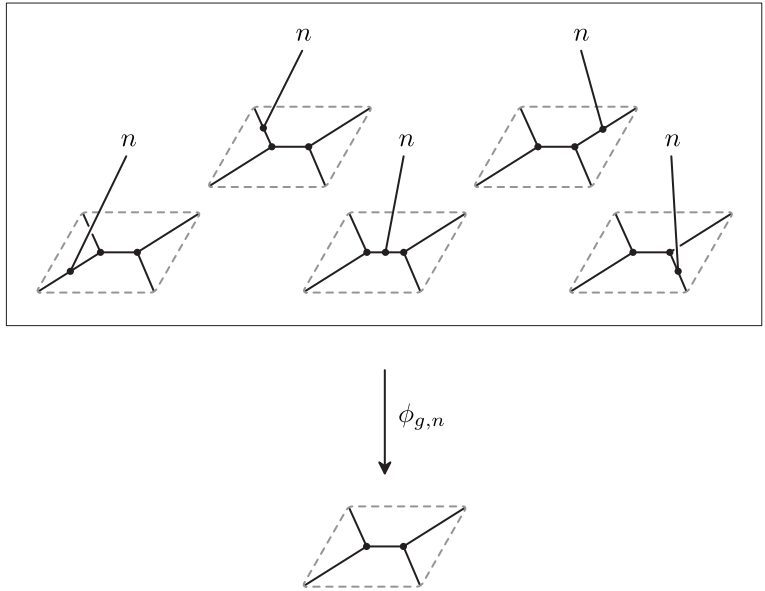}
\end{center}
\caption{The fiber over a vertex.}
\label{phi_inv_vert}
\end{figure}

Notice that two vertices of $\phi_{g,n}^{-1}(\Gamma)$ span an edge in 
$\phi_{g,n}^{-1}(\Gamma)$ if and only if they are graphs obtained from 
$\Gamma$ by inserting the new edge on adjacent edges of $\Gamma$, i.e. 
if and only if they can be obtained from one another by sliding the new
trivalent vertex, from one edge of $\Gamma$ to the adjacent one, through
their common trivalent vertex. As explained in Figure~\ref{edge_sliding},
this sliding is in fact an $F$ move performed on the marked edge. Then the
connectedness of $\phi_{g,n}^{-1}(\Gamma)$ immediately follows from that 
of $\Gamma$.

\begin{figure}[htb]
\begin{center}
\epsfig{file=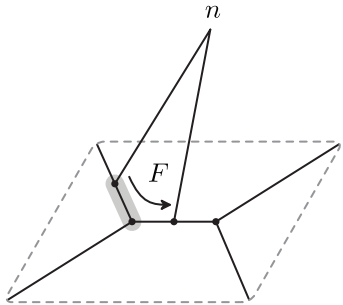}
\end{center}
\caption{The sliding of a trivalent vertex from one edge to the
  adjacent one is an $F$ move.}
\label{edge_sliding}
\end{figure}

The above observation allows us to define a canonical projection from
the 1-skeleton of $\phi_{g,n}^{-1}(\Gamma)$ to the graph $\Gamma$, once
they are barycentrically subdivided. Namely, we project any vertex
$\tilde\Gamma$ of $\phi_{g,n}^{-1}(\Gamma)$ to the barycenter of the edge
of $\Gamma$ split to get $\tilde\Gamma$ as a graph, and the barycenter
of any edge $F:\tilde\Gamma_1 \to \tilde\Gamma_2$ of
$\phi_{g,n}^{-1}(\Gamma)$ to the vertex of $\Gamma$ shared by the to
edges of $\Gamma$ split to obtain respectively $\tilde\Gamma_1$ and
$\tilde\Gamma_2$. Then, we extend the projection to a cellular map
between the two barycentric subdivisions in the obvious way.

Such projection induces a natural one-to-one correspondence between
all the paths in $\phi_{g,n}^{-1}(\Gamma)$ and those paths in the
barycentric subdivision of $\Gamma$ whose both ends are barycenters
of edges of $\Gamma$. Moreover, this correspondence respects composition
and sends loops to loops.

We need to show that any loop $\gamma$ in any fiber
$\phi_{g,n}^{-1}(\Gamma)$, with $\Gamma$ a vertex of $\S_{g,n-1}$,
can be contracted in $\S_{g,n}$ (by getting out of the fiber if
needed). The proof is by induction on the length $r$ of $\gamma$.

The base of the induction is the case $r=1$. In this case the loop
$\gamma$ corresponds to an $F$ move which, applied to a graph,
produces a new graph equivalent to the original one. This happens
when the $n$-th vertex is attached onto two adjacent edges $e_1$ and
$e_2$ of $\Gamma$, possibly coinciding, which are equivalent by an
automorphism of $\Gamma$. Figure~\ref{loop1} shows how it is possible 
to contract such a loop, exploiting the fact that it is one of the 
three edges of a triangle, whose other two edges coincide.

\begin{figure}[htb]
\begin{center}
\epsfig{file=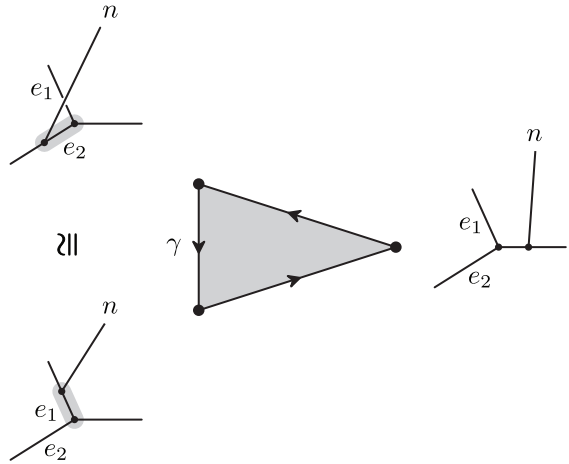}
\end{center}
\caption{A loop $\gamma$ of length $1$ in $\phi_{g,n}^{-1}(\Gamma)$ is
  contractible in $\S_{g,n}$.}
\label{loop1}
\end{figure}

For $r > 1$, the inductive step consists in proving that loop $\gamma$
is contractible in $\S_{g,n}$, assuming that the same holds for all the
loops of length $r-1$ in any fiber $\phi_{g,n}^{-1}(\Gamma)$.

\begin{figure}[b]
\begin{center}
\epsfig{file=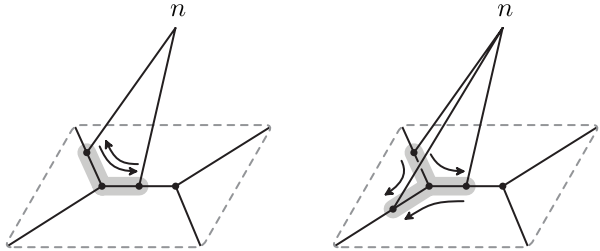}
\end{center}
\caption{Sequences of slidings inducing ``retraces'' in the loop
  $\hat{\gamma}$.}
\label{retracing}
\end{figure}

By means of the one-to-one correspondence introduced above, we can
interpret $\gamma$ as a loop of edges $\hat\gamma$ in $\Gamma$.

The loop $\hat\gamma$ may retrace one of its edges, i.e. may contain
an edge of the barycentric subdivision of $\Gamma$ followed by the
same edge with the opposite orientation. This may happen in one of
the two cases depicted in Figure~\ref{retracing}. In the first
situation, the retrace appears in $\gamma$ as well, and may be
canceled. The second situation represents a triangle in
$\phi_{g,n}^{-1}(\Gamma)$. By homotoping over the corresponding
$2$-cell, the retrace may be canceled.

After canceling all the retraces, we can assume that the length of
$\gamma$ and $\hat\gamma$ coincide.

Now, let us suppose that $\hat\gamma$ is not injective: this means
that an edge of $\Gamma$ appears two times in the loop $\hat\gamma$.
If this is the case, then the loop $\hat\gamma$ can be decomposed
into two loops of smaller length, inducing a similar decomposition
on $\gamma$. By induction on the length we are done.

Finally, we are reduced to the case when $\hat\gamma$ is an injective 
loop, i.e. it is a circuit in $\Gamma$ (see Figure~\ref{circuit}).

\begin{figure}[htb]
\begin{center}
\epsfig{file=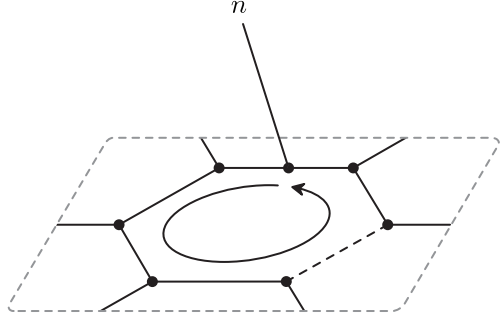}
\end{center}
\caption{A circuit of length $r$ in $\Gamma$.}
\label{circuit}
\end{figure}

We choose a vertex $\tilde\Gamma$ of $\gamma$. $\tilde\Gamma$ is a vertex
of $\S_{g,n}$, obtained from $\Gamma$ inserting a new edge with a free
end in the point $n$ and the other end creating a new trivalent vertex
on the edge $e\in\hat{\gamma}$. All the remaining vertices of $\gamma$
are obtained from $\Gamma$ inserting the new edge on the $r-1$ edges of
$\hat{\gamma}-\{e\}$. Let us denote by $\tilde\Gamma_1$ the vertex which
precedes $\tilde\Gamma$ and by $\tilde\Gamma_2$ the one following it
(with respect to the orientation of $\gamma$). Those vertices are
obtained from $\Gamma$ by attaching the new edge on the two edges of
$\hat{\gamma}$ adjacent to edge $e$ (see Figure~\ref{loop_reduction1}).

\begin{figure}[htb]
\begin{center}
\epsfig{file=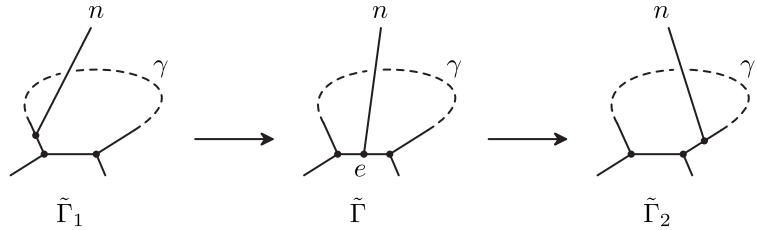}
\end{center}
\caption{Three consecutives vertices of a path $\gamma \in
  \phi_{g,n}^{-1}(\Gamma)$.}
\label{loop_reduction1}
\end{figure}

At this point, we apply to all vertices in $\gamma-\{\Gamma\}$ the
$F$ move performed on the edge $e$ ($F_e$ in the following), getting
from each of them a graph belonging to $\phi_{g,n}^{-1}(F_e(\Gamma))$. 
The graph $F_e(\Gamma)$ has a circuit of length $r-1$, and the vertices 
obtained above are those of the corresponding loop of $\phi_{g,n}^{-1}
(F_e(\Gamma))$ ($\gamma_e$ in Figure~\ref{loop_reduction2}).

The $F$ move which connects any pair of adjacent vertices in $\gamma
- \{\Gamma\}$, together with the two moves $F_e$ departing from such
vertices and with the $F$ move connecting their images in
$\phi_{g,n}^{-1}(F_e(\Gamma))$, bounds a DC square (see
Figure~\ref{loop_reduction2}).

\begin{figure}[htb]
\begin{center}
\epsfig{file=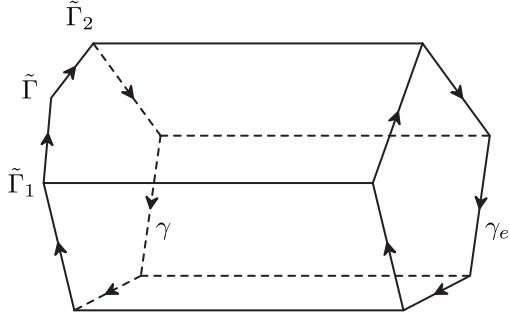}
\end{center}
\caption{Homotoping a loop of length $r$ to a loop of length $r-1$.
} \label{loop_reduction2}
\end{figure}

\break

Let us consider the pentagonal loop with vertices $\tilde\Gamma_1$,
$\tilde\Gamma$, $\tilde\Gamma_2$ (all belonging to $\phi_{g,n}^{-1}
(\Gamma)$), $F_e(\tilde\Gamma_1)$ and $F_e(\tilde\Gamma_2)$ (in 
$\phi_{g,n}^{-1}(F_e(\Gamma))$). As illustrated in 
Figure~\ref{loop_reduction3}, such loop can be filled in with
triangles and pentagons.

\begin{figure}[htb]
\begin{center}
\epsfig{file=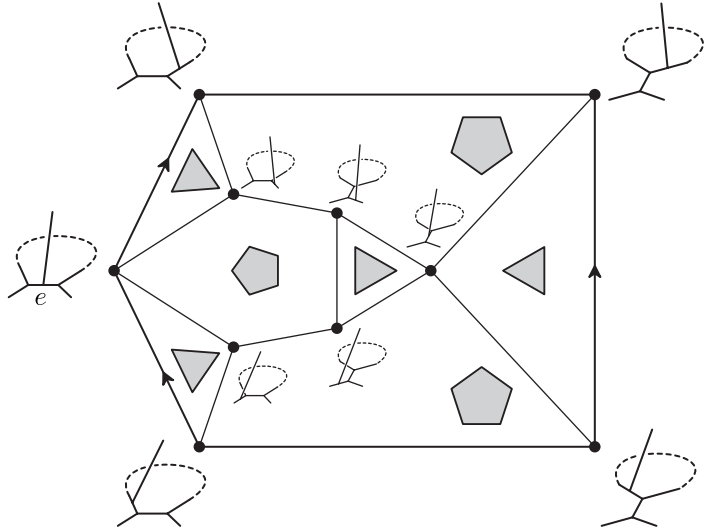}
\end{center}
\caption{The pentagonal loop may be filled in by triangles and pentagons.}
\label{loop_reduction3}
\end{figure}

The $r$-sided loop can therefore be homotoped to the $(r-1)$-sided one, 
which is in turn contractible by the inductive hypothesis, from which 
the thesis follows.

\smallskip
\noindent  {\bf\em Claim 2.} {\em The condition on the lifting of
edges is satisfied.}
\smallskip

Let $E:\Gamma \to \Gamma'$ be an edge of $\S_{g,n-1}$. Such $E$ is
an $F$ move performed on an edge $e$ (with distinct endpoints) of
$\Gamma$: hence, $E$ has $3g+2n-6$ liftings to $\S_{g,n}$, as many
as the number of edges in $\Gamma-\{e\}$. Let $E_1:\Gamma_1 \to
\Gamma'_1$ and $E_2:\Gamma_2 \to \Gamma'_2$ be two of these
liftings. As $\phi_{g,n}^{-1}(\Gamma)$ is connected, there is a
simple path $\gamma$ in $\phi_{g,n}^{-1}(\Gamma)$ connecting
$\Gamma_1$ to $\Gamma_2$. The graph $\Gamma_1$ is obtained from
$\Gamma$ by inserting a new trivalent vertex on an edge $f$ which is
different from $e$, and connecting it with the new free end $n$. The
path $\gamma$ consists of consecutive $F$ moves in the fiber, i.e.
consecutive slidings of the new trivalent vertex. We concentrate on
the first edge of $\gamma$, the one departing from $\Gamma_1$. Such
edge represents a sliding of the new vertex from edge $f$ to edge
$g$, through the vertex they have in common (i.e. it is an $F$ move
in the fiber $\phi_{g,n}^{-1}(\Gamma)$): denoted by $\Gamma_3$ the
resulting graph, still belonging to $\phi_{g,n}^{-1}(\Gamma)$, three
possibilities are given:
\begin{itemize}
\item[(i)]  $f$ and $e$ are disjoint;
\item[(ii)] $f$ and $e$ intersect at a vertex, and $g\neq e$;
\item[(iii)]$f$ and $e$ intersect at a vertex, and $g=e$.
\end{itemize}

In the first case, $E$ admits a lifting $E_3:\Gamma_3 \to \Gamma'_3$,
and $\Gamma'_3$ is obtained from $\Gamma'_1$ with the same sliding used
to get $\Gamma_3$ from $\Gamma_1$. The two liftings of $E$ and the two
slidings bound a DC square, and the resulting loop is therefore
contractible.

In the second case, $E$ still admits a lifting $E_3:\Gamma_3 \to
\Gamma'_3$, but to get $\Gamma'_3$ from $\Gamma'_1$ we need to perform
two slidings along consecutive vertices. Then, the two liftings of $E$
and the three slidings bound a pentagon, and the resulting loop is once
again contractible (see Figure~\ref{edge_lifting1}).

\begin{figure}[htb]
\begin{center}
\epsfig{file=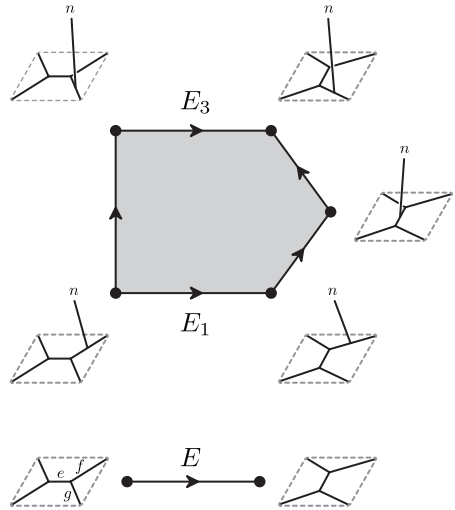}
\end{center}
\caption{Lifting of an edge, Case (ii).}
\label{edge_lifting1}
\end{figure}

\begin{figure}[htb]
\begin{center}
\epsfig{file=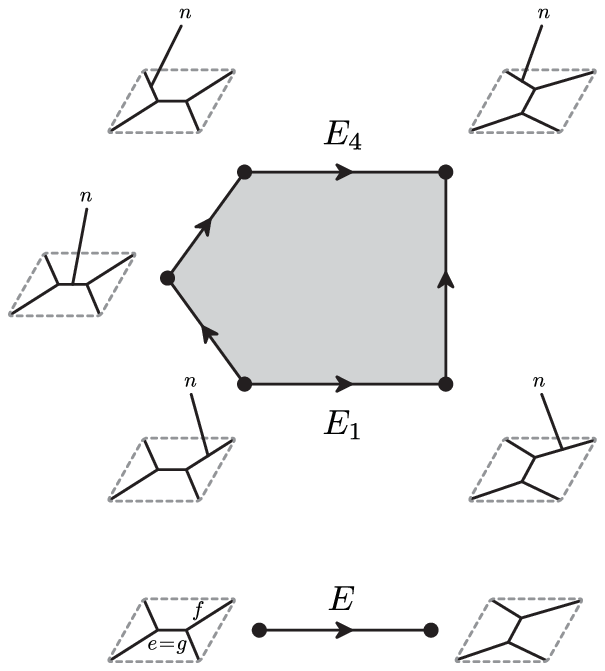}
\end{center}
\caption{Lifting of an edge, Case (iii-a).}
\label{edge_lifting2}
\end{figure}

\begin{figure}[htb]
\begin{center}
\epsfig{file=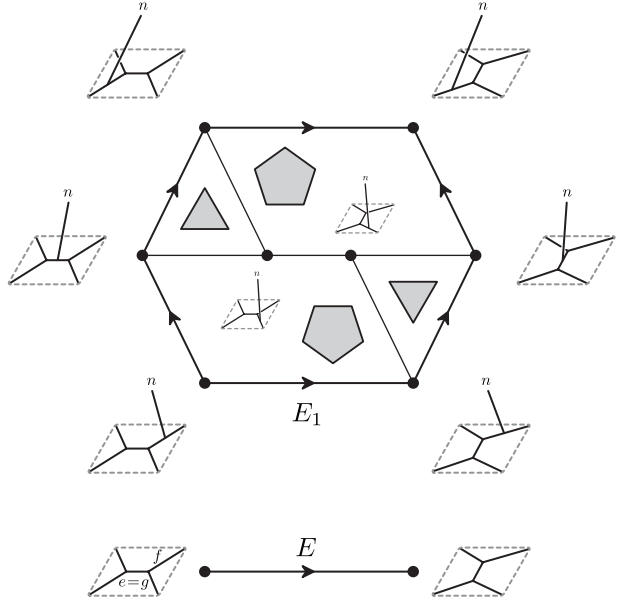}
\end{center}
\caption{Lifting of an edge, Case (iii-b).}
\label{edge_lifting3}
\end{figure}

In the third case, instead, $\Gamma_3$ does not support the $F$ move
that lifts $E$. However, if we slide further, by walking through the
following edge in the path $\gamma$, we end up with a graph $\Gamma_4$ 
that supports the move. Denoted by $E_4:\Gamma_4 \to \Gamma'_4$ the 
lifting of $E$ starting in $\Gamma_4$, two possibilities are given: 
either the graph $\Gamma'_4$ turns out to be related to $\Gamma'_1$ 
by a single sliding, or two slidings are needed to get $\Gamma'_1$ 
from $\Gamma'_4$. In the former case, the two liftings of $E$ and the 
three slidings still bound a pentagon, resulting in a contractible 
loop (see Figure~\ref{edge_lifting2}). In the latter case, the
two liftings and the four slidings give rise to an hexagon.
Nevertheless, such an hexagon may be subdivided into triangles and
pentagons, as shown in Figure~\ref{edge_lifting3},
being therefore contractible.

The same argument applies to all consecutive edges of the path $\gamma$.
We get this way a path $\gamma'$ in $\phi_{g,n}^{-1}(\Gamma')$ connecting
$\Gamma'_2$ to $\Gamma'_1$, and the loop bounded by $\gamma$, $\gamma'$,
$E_1$ and $E_2$ is contractible by construction. Hence Claim 2 is proved.

\smallskip

This concludes the proof of Proposition~\ref{RIDUZIONE}.
\end{proof}

\subsection{Dealing with the closed case}
\label{casochiuso}\strut\nobreak\medskip

A further step needed to prove Theorem~\ref{MMM}, is the definition of a
map $$\psi_g:\S_{g-1,2} \to \S_{g,0}\;.$$

\begin{figure}[htb]
\begin{center}
\epsfig{file=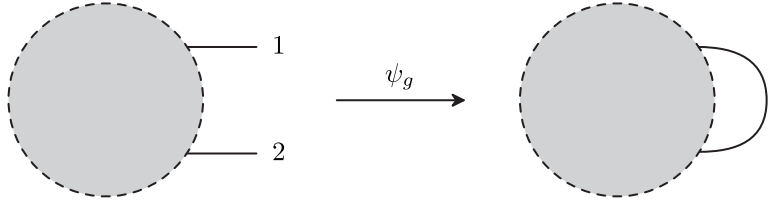}
\end{center}
\caption{Definition of $\psi_g$ on $V(\S_{g-1,2})$.}
\label{psi_vert}
\end{figure}

We define this map as follows. If $\Gamma$ is a vertex of $\S_{g-1,2}$, 
then $\psi_{g}(\Gamma)$ is the graph obtained from $\Gamma$ by attaching 
its two free ends together, as shown in Figure~\ref{psi_vert}. Such a 
graph turns out to have $2g-2$ trivalent vertices and no free ends, i.e. 
it is a vertex of $\S_{g,0}$. Now, every edge of $\S_{g-1,2}$ is an $F$ 
move connecting two graphs $\Gamma_1,\Gamma_2\in V(\S_{g-1,2})$. Then
there exists an $F$ move between $\psi_{g}(\Gamma_1)$ and
$\psi_{g}(\Gamma_2)$, and we define $\psi_{g}(F)$ to be such a move.
 The definition of $\psi_{g}$ on the 2-cells is straightforward.

The map $\psi_g$ we have just defined is cellular. It allows us to obtain
the second inductive step for the proof of Theorem \ref{MMM}.

\begin{PROP}\label{CASO_CHIUSO}
If the complex $\S_{g-1,2}$ is simply connected, then also $\S_{g,0}$
is simply connected, for every $g > 1$.
\end{PROP}

\begin{proof}
We observe that $\psi_g$ can be easily seen to be surjective, by the
same argument used for the surjectivity of $\phi_{g,n}$, except that here
the free edge is the one to be cut to get liftings of cells. Then the
connectedness of $\S_{g-1,2}$ implies that of $\S_{g,0}$, and the
proposition immediately follows from the following claim.

\smallskip
\noindent {\bf\em Claim 1.} {\em For any loop of edges $\gamma$ in
$\S_{g,0}$ there exists a loop of edges $\tilde\gamma$ in $\S_{g-1,2}$,
such that $\psi_g \circ \tilde\gamma$ is homotopic to $\gamma$ in
$\S_{g,0}$.}
\smallskip

In order to prove this claim, we need another claim.

\smallskip
\noindent {\bf\em Claim 2.} {\em For any vertex $\Gamma$ of $\S_{g,0}$
and any two vertices $\Gamma_1$ and $\Gamma_2$ in $\psi_g^{-1}(\Gamma)$,
there exists a path of edges $\delta$ in $\S_{g-1,2}$ between $\Gamma_1$
and  $\Gamma_2$, such that the loop $\lambda = \psi_g \circ \delta$ is
homotopically trivial in $\S_{g,0}$.}
\smallskip

\begin{figure}[b]
\begin{center}
\epsfig{file=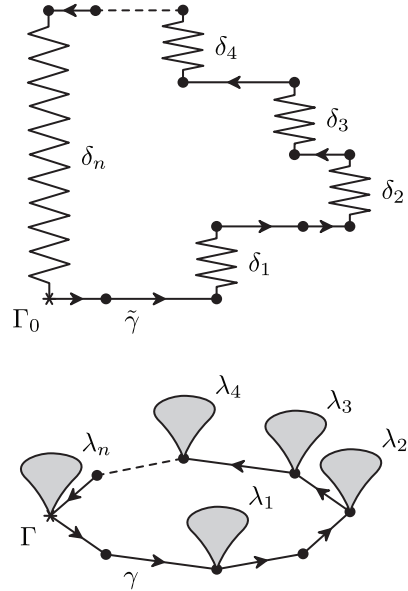}
\end{center}
\caption{The path $\tilde\gamma$ and its projection in $\S_{g,0}$.}
\label{contracting_loop_g0}
\end{figure}

We first prove Claim 1 assuming Claim 2.
We choose a vertex $\Gamma$ of $\S_{g,0}$ as the base point for the loop
$\gamma$, and denote by $\Gamma_0$ the vertex of $\psi_{g}^{-1}(\Gamma)$,
which is obtained from $\Gamma$ by cutting it along an edge $e_0$.

We choose $\Gamma_0$ as the basepoint of the lifting $\tilde\gamma$
and consider the first edge of $\gamma$, i.e. the one departing from
$\Gamma$. Such edge is an $F$ move performed on an edge $e$ of
$\Gamma$. If $e \neq e_0$, then we can lift the move $F_e$ to a move
on $\Gamma_0$. On the contrary, if $e = e_0$ then the move $F_e$
cannot be lifted to a move on the graph $\Gamma_0$ and we apply
Claim 2 in order to connect $\Gamma_0$ with a different vertex of
$\psi_g^{-1}(\Gamma)$ on which the move $F_e$ can be lifted.

We go on to construct our lifting, by following $\gamma$ edge after
edge and iterating the procedure described for the first edge.

In this way we end up with a path joining $\Gamma_0$ with some other
vertex of $\psi_g^{-1}(\Gamma)$ and we can close it to get the desired
loop $\tilde\gamma$ by applying once again Claim 2.
The scenario is that of Figure~\ref{contracting_loop_g0}.

Then $\psi_g \circ \tilde\gamma$ differs from $\gamma$ only for the
insertion of some homotopically trivial loops $\lambda_i = \psi_g \circ
\delta_i$, one for each application of Claim 2. Therefore, $\psi_g \circ
\tilde\gamma$ is homotopic to $\gamma$ as required by Claim 1.

Now we pass to prove Claim 2. Let $\Gamma_1$ and $\Gamma_2$ be
obtained from $\Gamma$ respectively by cutting two different edges
$e_1$ and $e_2$.

If $e_1$ and $e_2$ share both the ends, then $\Gamma_1 = \Gamma_2$
(see rightmost graphs in Figure \ref{contracting_psi_delta}) and there is
nothing to prove.

\begin{figure}[htb]
\begin{center}
\epsfig{file=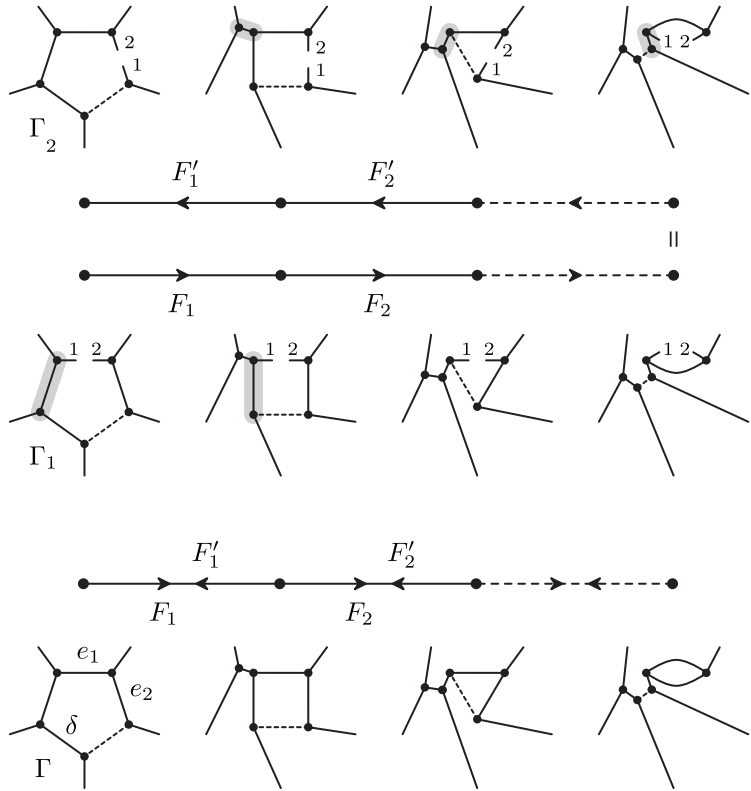}
\end{center}
\caption{The path $\psi_g \circ \tilde{\delta}$ is contractible in
  $\S_{g,0}$.}
\label{contracting_psi_delta}
\end{figure}

Otherwise, if $e_1$ and $e_2$ do not share both ends but are still
adjacent, then the fact that both can be cut without disconnecting
$\Gamma$ ensures the existence of a simple path of edges $\delta$ in
$\Gamma - \{e_1,e_2\}$ joining two different vertices of $e_1$ and
$e_2$, as in Figure \ref{contracting_psi_delta}. The same figure suggests
how to construct the desired path $\delta$ as a sequence of liftings
of $F$ moves performed in the order on the egdes of $\delta$,
followed by a sequence of liftings of their inverses in the reversed
order. Then the loop $\lambda = \psi_g \circ \delta$ turn out to be
homotopically trivial by construction.

Finally, if $e_1$ and $e_2$ are not adjacent, then we consider a
minimal path of edges $\epsilon$ in $\Gamma$ between a (trivalent)
vertex $v_1$ of $e_1$ and a vertex $v_2$ of $e_2$. By performing $F$
moves on the edges of $\epsilon$ in the order, we can drag $v_1$
until it coincides with $v_2$. In this way we get a path of edges
$\hat\epsilon$ in $\S_{g,0}$ connecting $\Gamma$ with a new graph
$\Gamma'$ where $e_1$ and $e_2$ are adjacent. Since neither $e_1$
nor $e_2$ can appear in $\epsilon$, we can lift $\hat\epsilon$ to
paths $\hat\epsilon_1$ and $\hat\epsilon_2$ in $\S_{g-1,2}$
respectively starting from $\Gamma_1$ and $\Gamma_2$. Now, calling
$\Gamma'_1$ and $\Gamma'_2$ the end points of these two paths in
$\psi_g^{-1}(\Gamma')$, we are reduced to the previous case.
Therefore, we can find a path of edges $\delta'$ between $\Gamma'_1$
and $\Gamma'_2$ such that the loop $\lambda' = \psi_g \circ \delta'$
is contractible in $\S_{g,0}$. Then we put $\delta = \hat\epsilon_1
\delta' \hat\epsilon_2^{-1}$ and observe that once again the loop
$\lambda = \psi_g \circ \delta = \epsilon \lambda' \epsilon^{-1}$ is
homotopically trivial by construction.
\end{proof}

\subsection{The decorated combinatorial structures}
\label{stilde}\strut\nobreak\medskip

As anticipated, we conclude this section with the construction of
the complex $\widetilde{\S}_{g,n}$, codifying all of {\em decorated
combinatorial structures}, i.e. the combinatorial structures of pant decompositions whose curves are ordered.

Let $V(\widetilde{\S}_{g,n})$ be the set of all connected graphs
with $n$ univalent vertices and $2g-2+n$ trivalent vertices,
equipped with an ordering of the $3g-3+n$ edges connecting two
trivalent vertices. For those graphs we may define a {\em decorated
$F$ move}, naturally lifting the combinatorial $F$ move defined in
Section \ref{CS}, simply requiring that the new edge created by the
$F$ move inherits its number by the old one. In order to emphasize
that the $F$ move is performed on the $i$-th edge, we denote it by
$F_i$. Moreover, if two decorated graphs are identical except for
the ordering of the edges, which differ by a transposition
$(ij)\in\Sigma_{3g-3+n}$, we connect them with an edge, labeled by
$\tau_{ij}$. We complete the construction by adding 2-cells of three
types: {\em combinatorial}, {\em algebraic} and {\em mixed} ones.

\begin{figure}[htb]
\begin{center}
\epsfig{file=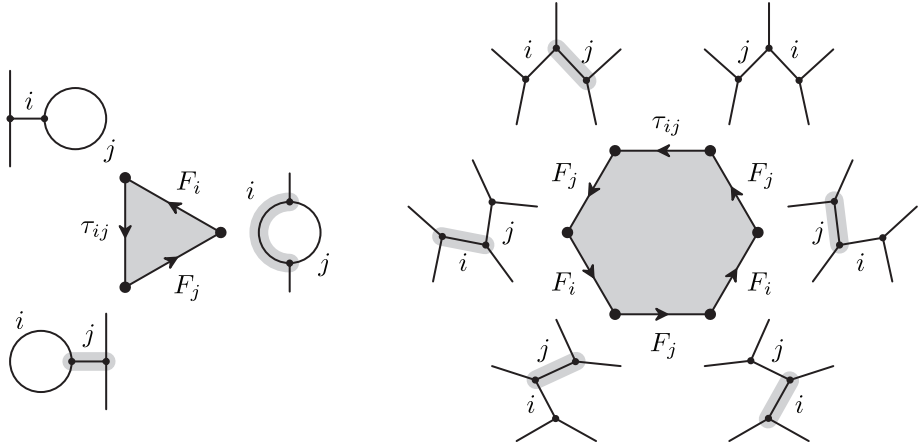}
\end{center}
\caption{A bigon and a pentagon of decorated $F$ moves.}
\label{decorated_bigon_pentagon}
\end{figure}

The combinatorial 2-cells translate in terms of decorated
combinatorial structures the bigons, triangles, squares and
pentagons of $\S_{g,n}$. Notice that the shape of the bigons and
that of the pentagons change when lifted to the decorated setting,
as depicted in Figure~\ref{decorated_bigon_pentagon}, while
all the other combinatorial cells maintain their shape.

The algebraic $2$-cells correspond to the relation in the group
$\Sigma_{3g-3+n}$, thus they are squares $\tau_{ij}\tau_{hk}\tau_{ij}
= \tau_{lm}$, for any couple of different transpositions $(ij)$ and 
$(hk)$ in $\Sigma_{3g-3+n}$, where $(lm)=(ij)(hk)(ij)$.

Finally, the mixed two cells are squares $F\tau_{ij}=\tau_{ij}F$,
telling that $F$ and $\tau$ moves commute.

The complex $\widetilde{\S}_{g,n}$ built like that has an obvious
projection on $\S_{g,n}$, which is defined on the vertices by
forgetting the ordering and extends to edges and 2-cells in a
natural way. Proving that such map satisfies all the hypothesis of
Proposition \ref{teoremazero} is straightforward, and we leave it to
the reader. The application of Proposition \ref{teoremazero} thus
ensures that $\widetilde{\S}_{g,n}$ is connected and simply
connected.

\section{The complex of pant decompositions}
\label{decomposizioni in pantaloni}

In this section we want to lift the complex $\S_{g,n}$ to a new
complex $\R_{g,n}$, whose vertices are in one-to-one correspondence
with the {\em pant decompositions} of the surface $\Sigma_{g,n}$.
This result is achieved in two steps, described in the two
subsections: first of all, we define a complex
$\widetilde{\R}_{g,n}$, whose vertices are in one-to-one
correspondence with the decorated pant decompositions of
$\Sigma_{g,n}$, in such a way that the natural action of the mapping
class group $\M_{g,n}$ on the decorated pant decompositions extends
to an action on $\widetilde{\R}_{g,n}$ and the quotient
$\widetilde{\R}_{g,n}/\M_{g,n}$ coincides with
$\widetilde{\S}_{g,n}$. Provided we take care that the natural
projection $r_{g,n}:\widetilde{\R}_{g,n} \to \widetilde{\S}_{g,n}$ 
satisfies the conditions of Proposition~\ref{teoremazero}, the new
complex turns out to be simply connected. Finally, we exploit the
natural projection $t_{g,n}:\widetilde{\R}_{g,n} \to \R_{g,n}$ to
prove that $\R_{g,n}$ is simply connected as well.

\subsection{The decorated pant decompositions}
\label{decorated_pant_decompositions}\strut\nobreak\medskip

As anticipated in the introduction, we need now a presentation for
all the mapping class groups. In order to produce such presentation,
we use the results of~\cite{Benv}. In such paper a machinery for
finding presentations is built; its input is a presentation for the
mapping class group of the sphere with $4$ and $5$ boundary
components, and for the torus with $1$ and $2$ boundary components
({\em sporadic surfaces}). Moreover, \cite{Benv} also shows that
such a machinery produces a presentation of any known ``style'',
provided the input is chosen according to the same ``style''. In
order to perform the explicit calculations shown in this section, it
is convenient to use a presentation in terms of Dehn twists,
described by Gervais in~\cite{Gervais0}, that we may obtain with the
above recalled method starting from the {\em Dehn twist style}
presentations for the sporadic surfaces.

The generators in such a presentation are the Dehn twists along all
simple closed curves in $\Sigma_{g,n}$, while the relations belong
to three simple types, {\em braids}, {\em lanterns} and {\em
chains}.

Namely, we call {\em braids} the relations of the form
$$T_c=T_bT_aT_b^{-1},$$ where the curves $a$ and $b$ are such that $|a
\cap b| = 0,1$ or $|a \cap b| = 2_0$ (meaning they intersect in two
points with algebraic intersection zero), and $c=T_b(a)$.

\begin{figure}[htb]
\epsfig{file=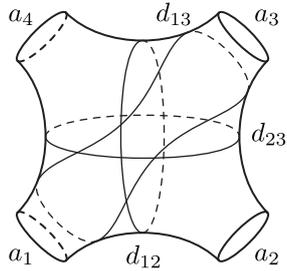} 
\caption{Lanterns are supported in subsurfaces homeomorphic to 
  $\Sigma_{0,4}$.}
\label{lantern}
\end{figure}

We call {\em lanterns} the relations like
$$T_{a_1}T_{a_2}T_{a_3}T_{a_4} = T_{d_{12}}T_{d_{23}}T_{d_{13}},$$
where the curves $a_i,\, d_{ij}$ are represented in
Figure~\ref{lantern}.

\begin{figure}[htb]
\epsfig{file=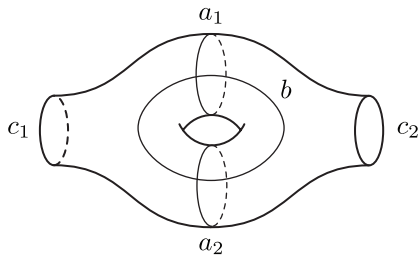}
\caption{Chains are supported in subsurfaces homeomorphic to 
  $\Sigma_{1,2}$.}
\label{chain}
\end{figure}

Finally, we name {\em chains} the relations
$$(T_{a_1}T_{b}T_{a_2})^4=T_{c_1}T_{c_2},$$
where $a_i,\, c_{j}$ are the curves depicted in Figure~\ref{chain}.

The vertices of $\widetilde{\R}_{g,n}$ have to be in a one to one
correspondence with the (infinitely many) decorated pant
decompositions of $\Sigma_{g,n}$, i.e.
$$V(\widetilde{\R}_{g,n}) = \{\text{decorated pant decompositions of
$\Sigma_{g,n}$}\}/\text{isotopy}$$

We now define a transformation between decorated pant decompositions, 
that we call again {\em $F$ move}. Let $D$ be the decorated pant 
decomposition given by an ordered family $\{a_1,\dots,a_{3g-3+n}\}$ of 
curves (up to isotopy). Let us consider one of these curves, $a_i$, and 
let $a'_i$ be any curve on $\Sigma$ such that $|a_i \cap a'_i| = 2_0$ and
$D'=\{a_1,\dots,a'_i,\dots, a_{3g-3+n}\}$ is still a decorated pant
decomposition. The move $F$ is defined as the transformation
$$D=\{a_1,\dots,a_i,\dots,a_{3g-3+n}\} \stackrel{F}{\longto}
D' =\{a_1, \dots,a'_i,\dots, a_{3g-3+n}\}.$$ In order to emphasize
that such $F$ move is performed on the $i$-th curve of $D$, we denote 
it by $F_i$.

\begin{figure}[htb]
\begin{center}
\epsfig{file=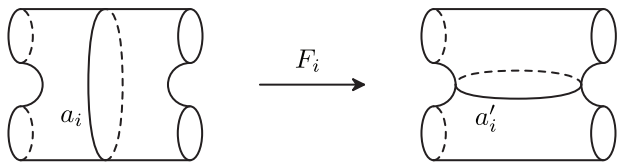}
\end{center}
\caption{The $F$ move between decorated pant decompositions.}
\label{Fmove_decomp}
\end{figure}

The move, depicted in Figure~\ref{Fmove_decomp}, is a transformation 
between decorated pant decompositions of $\Sigma_{g,n}$, supported in a
subsurface homeomorphic to a sphere with four boundary components, i.e. 
it is a local move.

We connect two vertices in $V(\widetilde{\R}_{g,n})$ by an edge $F$ if 
the corresponding decompositions are related to one another by an $F$ 
move. Moreover, we insert an edge $\tau_{ij}$ between two vertices if 
the corresponding decompositions are identical but for the ordering of 
the curves, which differ by the transposition $\tau_{ij}, \quad i,j=1,
\ldots,3g-3+n, i\neq j$. Finally, we insert an edge $T_a$ between two 
vertices if the corresponding decompositions are related to one another 
by the Dehn twist along a simple closed curve $a$.

To the 1-dimensional complex obtained above, we now add 2-cells of
{\em combinatorial}, {\em algebraic}, {\em topological} and {\em
mixed} type.

\begin{figure}[b]
\begin{center}
\epsfig{file=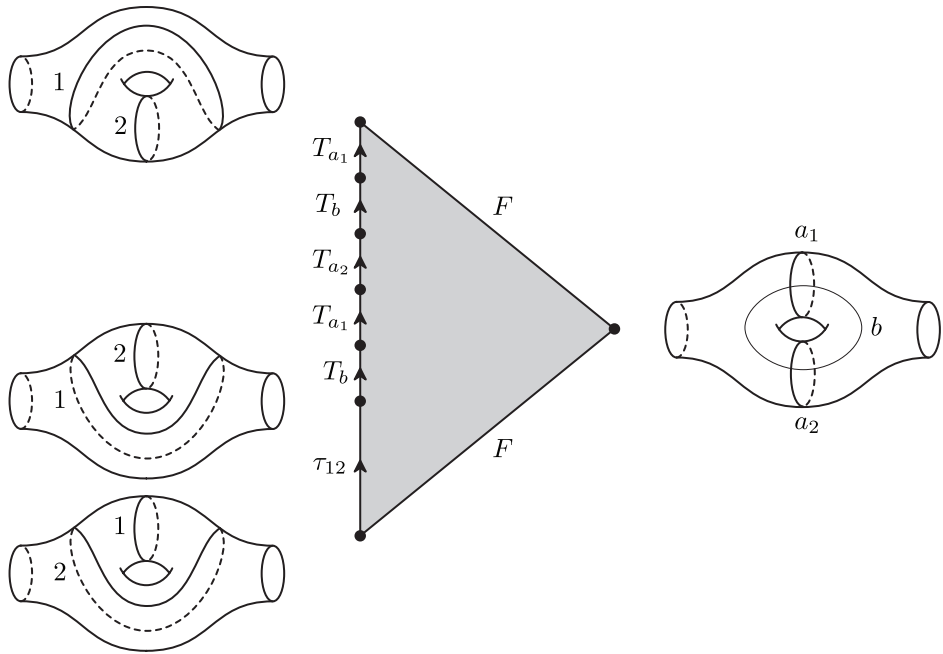}
\end{center}
\caption{A bigon of decorated pant decompositions.}
\label{decorated_bigon_decomp}
\end{figure}

The combinatorial 2-cells rephrase in terms of pant decompositions the 
bigons, triangles, squares and pentagons of $\widetilde{\S}_{g,n}$. 
Once again, the shape of the bigons changes when lifted from the 
combinatorial setting to the new one, as depicted in 
Figure~\ref{decorated_bigon_decomp}, while all the other combinatorial 
cells maintain their shape.

The algebraic 2-cells correspond to the relations of the symmetric
group over $3g-3+n$ elements, hence they are squares of $\tau$ moves.

\begin{figure}[htb]
\begin{center}
\epsfig{file=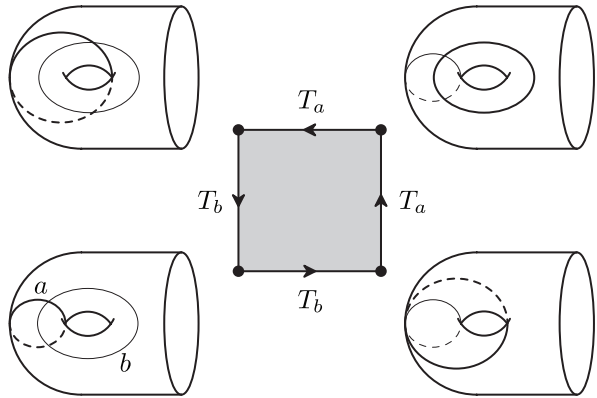}
\end{center}
\caption{A square of $T$ moves.} 
\label{hatcher}
\end{figure}

The topological 2-cells are the ones carried by the relations in the
Dehn twist presentation of $\M_{g,n}$ (braids, lanterns and chains),
together with the squares of $T$ moves shown in Figure~\ref{hatcher}. 
We also see that a twist, when performed on a curve which does not 
intersect any of the curves belonging to $D$, produces a loop of lenght 
$1$ based in the vertex of $\widetilde{\R}_{g,n}$ corresponding to $D$. 
If we then fill also these loops with 2-cells, we obtain a third type of
topological 2-cells.

\begin{figure}[htb]
\begin{center}
\epsfig{file=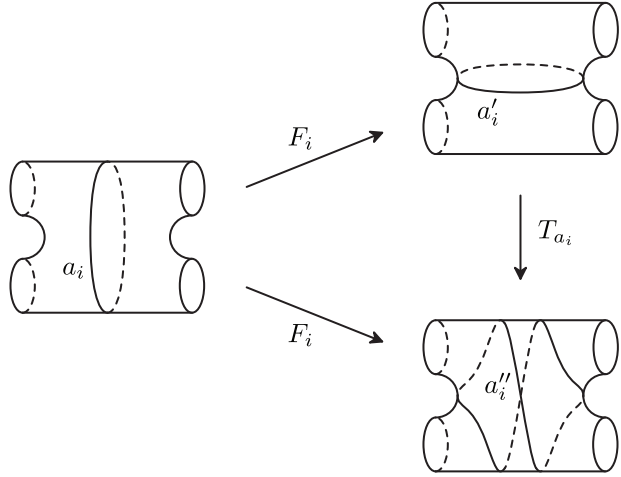}
\end{center}
\caption{Mixed triangles.} 
\label{mixed_triangle_decomp}
\end{figure}

Finally, the mixed 2-cells are the triangles shown in
Figure~\ref{mixed_triangle_decomp}, the squares telling that the $F$ 
moves commute with the Dehn twists and with the $\tau$ moves, and those 
telling that Dehn twist and $\tau$ moves commute as well.

We define $\widetilde{\R}_{g,n}$ to be the complex having
$V(\widetilde{\R}_{g,n})$ as the set of vertices, edges and 2-cells
as above. We may now state the following result, whose proof takes up 
the remainder of this subsection.

\begin{Th}\label{complesso_decorato}
The complex $\widetilde{\R}_{g,n}$ is simply connected.
\end{Th}

\begin{proof}
To prove the theorem, we consider the map
$r_{g,n}:\widetilde{\R}_{g,n} \to \widetilde{\S}_{g,n}$ defined as
follows. On the vertices, $r_{g,n}$ is the natural projection
associating to any decorated pant decomposition the corresponding
decorated combinatorial structure. As far as the edges are concerned, 
$r_{g,n}$ sends each $F$ and each $\tau_{ij}$ between two decompositions 
in the combinatorial $F$ and $\tau_{ij}$ between the corresponding
combinatorial structures, and contracts each Dehn twist $T_a$ to a
point. The map so defined extends to a map $r_{g,n}:\widetilde{\R}_{g,n} 
\to \widetilde{\S}_{g,n}$, which is surjective thanks to the combinatorial 
and the algebraic 2-cells and to the squares $F\tau=\tau F$ inserted in 
$\widetilde{\R}_{g,n}$.

We are then left to check that the remaining hypotheses of
Proposition~\ref{teoremazero} are fulfilled.

\smallskip
\noindent {\bf\em Claim 1.} {\em The fiber over each vertex is
connected and simply connected in $\T_{g,n}$.}
\smallskip

Let $\Gamma$ be a vertex of $\widetilde{\S_{g,n}}$, i.e a
combinatorial structure of decorated pant decomposition. Then,
$$r_{g,n}^{-1}(\Gamma) = \{\text{decorated pant decompositions with
combinatorial structure $\Gamma$}\}$$

Given one element of $r_{g,n}^{-1}(\Gamma)$, any other is obtained
from the selected one via the action of the mapping class group. The
set $r_{g,n}^{-1}(\Gamma)$ is therefore connected, since the edges
of $\widetilde{\R}_{g,n}$ include all the Dehn twists.

We need to prove that $r_{g,n}^{-1}(\Gamma)$ is simply connected.
Since every element $\gamma$ in $\M_{g,n}$ defines a corresponding
path in $r_{g,n}^{-1}(\Gamma)$, which is unique up to homotopy in
the fiber, such path will be denoted by $\gamma$ as well. Let
$\gamma$ be a loop in $r_{g,n}^{-1}(\Gamma)$, and let $D$ be its
basepoint. There are two possibilities: either $\gamma$ is a trivial
element in $\M_{g,n}$, or it is a nontrivial element of
$Stab(D)\subset\M_{g,n}$. In the former case, $\gamma$ is
contractible due to the topological 2-cells, representing the
relations of the mapping class group. In fact, $\gamma$ is a
relation in $\M_{g,n}$, then it follows from the braids, lanterns
and chains. Let then $\gamma$ represent a nontrivial element of
$Stab(D)\subset\M_{g,n}$. In particular, $\gamma$ sends each curve
of $D$ in itself, possibly inverting its orientation.

If $\gamma$ respects the orientation of all the curves in $D$, then
it may be expressed (i.e. is equivalent in $\M_{g,n}$) as the
product of Dehn twists performed over the curves $a_i$ of the
decomposition $D$. By means of the relations of $\M_{g,n}$, the loop
$\gamma$ may be homotoped in $r_{g,n}^{-1}(\Gamma)$ to the product
of the corresponding loops, i.e.
$$\gamma \sim \prod_i T_{a_i}.$$
Such loops are in turn contractible, because we cupped them off by
the corresponding 2-cells, thus ensuring that $\gamma$ is
contractible as well.

Let us suppose that $\gamma$ respects the ordering of the curves in
$D$, but it changes the orientation of one of them, $a_i$. Thus,
$a_i$ cannot be a separating curve of $\Sigma_{g,n}$: in fact, if
$a_i$ was separating, then $\gamma$ should switch the two connected
components of $\Sigma_{g,n}-a_i$. Then, $\gamma$ would switch at
least two curves in $D$, or two components of
$\partial\Sigma_{g,n}$, which is impossible.

Hence, $a_i$ is non-separating and two possibilities are given:
either $a_i$ bounds on both sides the same pant $P$, or $a_i$ bounds
a pant $P$ on one side and a pant $P'$ on the other.

In the first situation, the homeomorphism $\gamma$ may be
represented by the product $\omega_i\cdot \prod_j T_{a_j}$, where
the $a_j$'s are curves of $D$ and $\omega_i$ is the semitwist of $P$
relative to $a_i$. It is well known that such a semitwist may be
expressed by the product $T_bT_{a_i}^2T_b$, where $b$ is as shown in
Figure~\ref{hatcher}).
Hence the loop $\gamma$ is homotopic in $r_{g,n}^{-1}(\Gamma)$ to
the product of the corresponding loops, i.e.
$$\gamma \sim T_bT_{a_i}^2T_b\prod_j T_{a_j}.$$
The loop corresponding to the semitwist is contractible (indeed it
is the boundary of a square of $T$ moves), and so are the others,
ensuring this way that $\gamma$ is contractible as well.

In the second situation, $\gamma$ switches the two pants bounded by
$a_i$, $P$ and $P'$. Since $\gamma$ does not permute the curves of
$D$, also the remaining boundary components of $P$ and $P'$ must be
in common. Hence, $\Sigma_{g,n}=\Sigma_{2,0}$, $D$ is as shown in
Figure~\ref{rotation_20} and $\gamma$ is the
rotation of $\pi$ radians around the horizontal axis.

\begin{figure}[htb]
\begin{center}
\epsfig{file=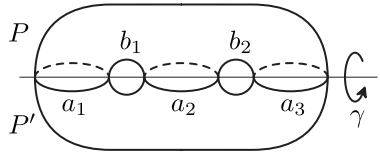}
\end{center}
\caption{The rotation $\gamma$.}
\label{rotation_20}
\end{figure}

We remark that all the loops we considered so far are contractible
without getting out of their own fiber, which is any one among the
$r_{g,n}^{-1}(\Gamma)$. Conversely, in the case of the rotation
depicted in Figure~\ref{rotation_20} it is required to get out of 
the fiber, as described in the following lemma.

\begin{LEMMA}
Let $\gamma$ be the homeomorphism of $\Sigma_{2,0}$ switching $P$
and $P'$ (notations as in Figure~\ref{rotation_20}). Then the 
corresponding loop of $r_{2,0}^{-1}(\Gamma)$ is contractible in
$\widetilde{\R}_{2,0}$.
\end{LEMMA}

\begin{proof}
The homeomorphism $\gamma$ may be expressed in terms of Dehn twists
as the product
$\gamma=(T_{b_2}T_{a_3}^{2}T_{b_2})^{-1}T_{b_1}T_{a_1}^{2}T_{b_1}$
(up to twists along the curves of $D$). Let us perform, starting
from $D$, an $F$ move along the curve $a_2$, and let us denote by
$D'$ the resulting pant decomposition. The homeomorphism $\gamma$,
applied to $D'$, gives rise to a contractible loop, as shown in
Figure~\ref{rotation_loop_20}.
\begin{figure}[htb]

\begin{center}
\epsfig{file=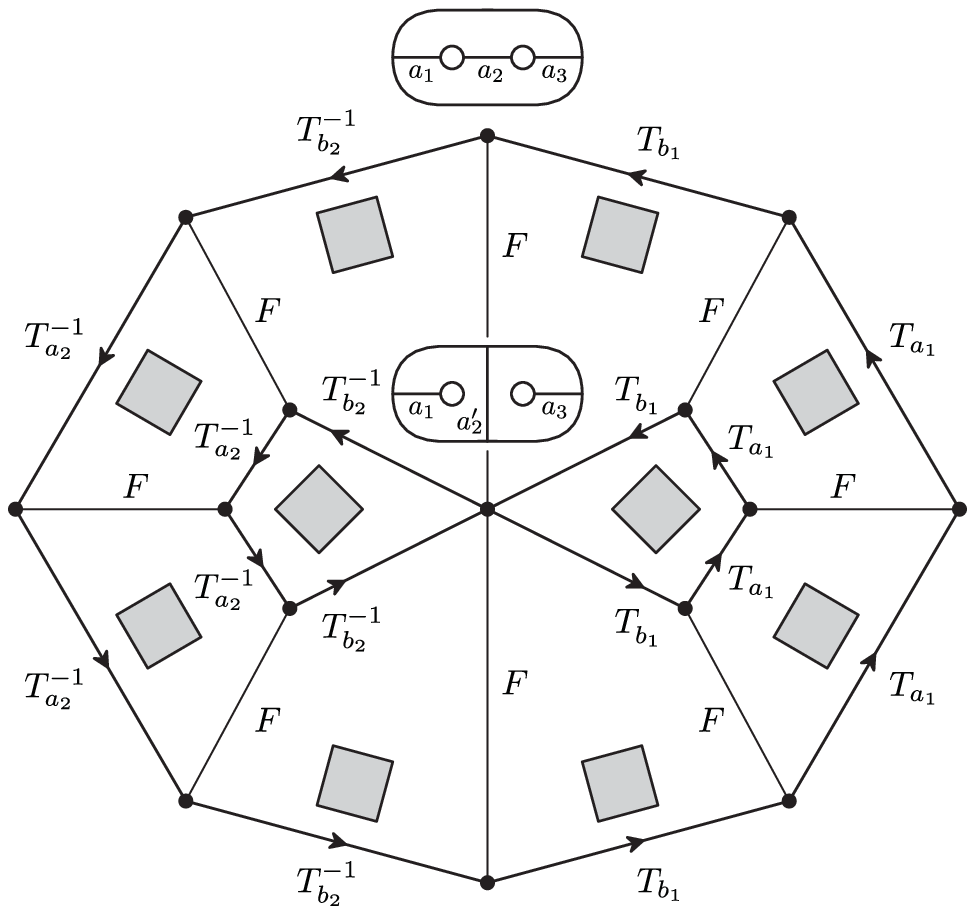}
\end{center}
\caption{Contracting the loop $\gamma$ in $\widetilde{\R}_{2,0}$.}
\label{rotation_loop_20}
\end{figure}

The same picture shows that the original loop $\gamma$ is homotopic,
by means of the mixed squares, to that contractible loop, thus
proving the thesis.

\end{proof}

In conclusion, if $\gamma$ inverts the orientation of more than one
curve of $D$, then either we are in the situation of the lemma
(hence $\gamma$ is contractible) or all curves whose orientation is
changed by $\gamma$ bound on both sides the same pant. In the latter
case, $\gamma$ may be expressed as the product of Dehn twists
performed over the curves $a_i$ and of semitwists $\omega_j$
relative to the $a_j$'s whose orientation is changed. By means of
the relations of $\M_{g,n}$, the loop $\gamma$ may be homotoped in
$r_{g,n}^{-1}(\Gamma)$ to the product of the corresponding loops,
i.e.
$$\gamma \sim \prod_i T_{a_i}\prod_j \omega_j.$$
Such loops are in turn contractible, since we cupped them off by the
corresponding 2-cells, and this concludes the proof of Claim 1.

\smallskip
\noindent  {\bf\em Claim 2.} {\em The condition on the lifting of
edges is satisfied.}
\smallskip

The edges of $\widetilde{\S}$ are $F$ and $\tau$ moves. Each of
these edges has infinitely many liftings to $\widetilde{\R}$, and
the lifting condition follows almost trivially by the mixed 2-cells.

Namely, let us consider a decorated combinatorial $F$ move, $F_i$,
performed on an edge $e_i$, that transforms the graph $\Gamma_1$
into the graph $\Gamma_2$. Let us consider two different liftings of
the $F_i$ move, $F'_i: D'_1 \to D'_2$ and $F''_i: D''_1 \to D''_2$.
Given a path $\delta_1$ in $r_{g,n}^{-1}(\Gamma_1)$ connecting
$D'_1$ to $D''_1$, the same path connects $D'_2$ to a $D_2$ in
$r_{g,n}^{-1}(\Gamma_2)$. Such $D_2$ coincides with $D''_2$ except,
possibly, for the $i-th$ curve. $D_2$ and $D''_2$ are then connected
by a power of $T_{a'_i}$. The path $\delta_2$, which is the
composition of $\delta_1$ with the suitable power of $T_{a'_i}$,
connects $D'_2$ to $D''_2$ in $r_{g,n}^{-1}(\Gamma_2)$, and the
square bounded by $F'_i, F''_i, \delta_1$ and $\delta_2$ is
contractible due to the mixed squares and triangles. If
$D'_1=D''_1$, then $D'_2$ and $D''_2$ simply differ by a multiple of
$T_{a'_i}$. Thus, the loop $F_i''T_{a'_i}^k\left(F_i'\right)^{-1}$
is contractible due to the mixed triangles.

The proof of the analogous condition for the lifting of edges of
type $\tau$ is straightforward, provided we remark that such moves
commute with all Dehn twists and that we inserted all the
corresponding DC squares.

\smallskip

All conditions of Proposition~\ref{teoremazero} are fulfilled. The
thesis of Theorem~\ref{complesso_decorato} has therefore been
demonstrated.

\end{proof}

\subsection{Back to pant decompositions}
\label{back}\strut\nobreak\medskip

The final step of our construction consists in the definition of a
complex $\R_{g,n}$, with all the pant decompositions of $\Sigma_{g,n}$ 
as the set of vertices, such that two maps $t_{g,n}:\widetilde{\R}_{g,n} 
\to \R_{g,n}$ and $p_{g,n}:\R_{g,n} \to \S_{g,n}$ are defined and the 
diagram 
$$\epsfig{file=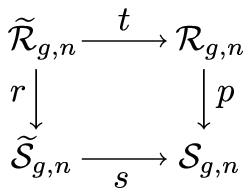}$$
(where $s$ is the natural projection) is commutative. The core
result of this last section is the proof that the simply
connectedness of $\R_{g,n}$ follows directly from that of
$\widetilde{\R}_{g,n}$.

\begin{figure}[htb]
\begin{center}
\epsfig{file=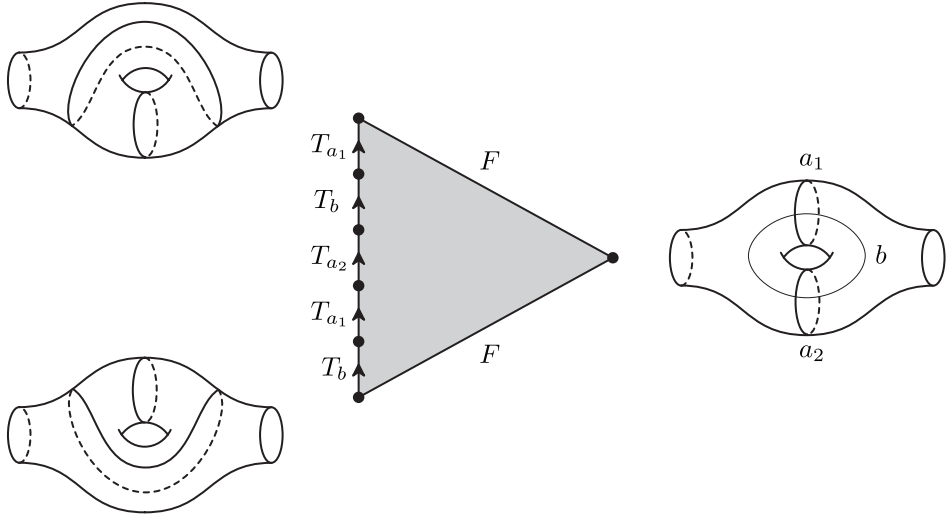}
\end{center}
\caption{A bigon of pant decompositions.}
\label{bigon_decomp}
\end{figure}

Let $V(\R_{g,n})$ be the set of all pant decompositions of 
$\Sigma_{g,n}$, considered up to isotopy. We connect two vertices in
$V(\R_{g,n})$ by an edge $F$ if the corresponding decompositions are
related to one another by an $F$ move (where the $F$ move is defined
as in $\widetilde{\R}_{g,n}$, just forgetting about the ordering of
the curves). Moreover, we insert an edge $T_a$ between two vertices
if the corresponding decompositions are related to one another by
the Dehn twist along a simple closed curve $a$. To the 1-dimensional
complex obtained above, we add 2-cells of {\em combinatorial} type
(bigons, triangles, squares and pentagons of $F$ moves), {\em
topological} type (braids, lanterns, chains, squares of $T$ moves
and one sided 2-cells corresponding to twists along the curves of a
decomposition) and {\em mixed} type (triangles as in
Figure~\ref{mixed_triangle_decomp}, squares $F T_a =T_a F$).
We remark that bigons and pentagons change their shape when passing
from the decorated setting to the non-decorated one (as depicted in
Figure~\ref{bigon_decomp}), while all the other cells
maintain the shape of the corresponding ones of
$\widetilde{\R}_{g,n}$.

We define $\R_{g,n}$ to be the complex having $V(\T_{g,n})$ as the
set of vertices, edges and 2-cells as above. The two maps
$t_{g,n}:\widetilde{\R}_{g,n} \to \R_{g,n}$ (which forgets about the
ordering of the curves) and $p_{g,n}:\R_{g,n} \to \S_{g,n}$ (which
forgets about the topological information, keeping track of the the
combinatorial one only) are well defined. Moreover, it is easy to
show that $p\circ t=s\circ r$. We may now state our core result.

\begin{Th}\label{complesso_intero}
The complex $\R_{g,n}$ is simply connected.
\end{Th}

\begin{proof}
Let $\gamma$ be a loop in $\R_{g,n}$, based at the decomposition
$D$. Hence $\gamma$ is a sequence of $F$ moves and Dehn twists,
transforming the decomposition $D$ into itself. We choose a
decorated pant decomposition $D_1\in t_{g,n}^{-1}(D)$ and we lift
the loop $\gamma$ edge after edge, starting form $D_1$. We get in
this way a path $\tilde{\gamma}$ in $\widetilde{\R}_{g,n}$, whose
second end is a decomposition $D_2$, still belonging to
$t_{g,n}^{-1}(D)$. The curves in $D_1$ and $D_2$ are identical,
possibly except for their enumeration. Hence there exists a path of
$\tau$ moves connecting $D_2$ to $D_1$ in $t_{g,n}^{-1}(D)$. Such
path allows us to close the path $\tilde{\gamma}$ to a loop, based
at $D_1$, which will be denoted by $\bar{\tilde{\gamma}}$. As
$\widetilde{\R}_{g,n}$ is simply connected, $\bar{\tilde{\gamma}}$
is contractible, i.e. it may be homotoped to a point by means of the
2-cells of $\widetilde{\R}_{g,n}$. Projecting the loop
$\bar{\tilde{\gamma}}$ and the homotopy in $\R_{g,n}$ with the map
$t$, we get that the loop $\gamma=t(\bar{\tilde{\gamma}})$ is
contractible as well, thus proving the Theorem.
\end{proof}

\begin{REM}
The topological cells provided by the Dehn style presentation are
supported in subsurfaces homeomorphic to sporadic surfaces, living
at the first and second Grothendieck floor (i.e. such that
$3g-3+n=1,2$). The same is true for all the other cells of
$\R_{g,n}$. Hence, $\R_{g,n}$ turns out to have a presentation with
generators and relations supported in sporadic surfaces,
illustrating in this way the simplified version of the Grothendieck
conjecture recalled in the introduction.
\end{REM}

\bibliographystyle{amsplain}

\end{document}